\documentclass[12pt]{amsart}
\usepackage{latexsym, amssymb, amsmath, amscd}
 \usepackage{eucal}

    \newtheorem{rema}{Remark}[section]
    \newtheorem{propo}[rema]{Proposition}

   \newtheorem{theo}[rema]{Theorem}
   \newtheorem{def-theo}[rema]{Definition-Theorem}
 \newtheorem{conj}[rema]{Conjecture}
   \newtheorem{defi}[rema]{Definition}
    \newtheorem{lemma}[rema]{Lemma}
    \newtheorem{corol}[rema]{Corollary}
     \newtheorem{exam}[rema]{Example}
  \newtheorem{rmk}[rema]{Remark}

	\newcommand{\nno}{\nonumber}

	\newcommand{\p}{\partial}

 \newcommand{\pf}{{\it Proof:}\hspace{2ex}}
 
 \newcommand{\epfv}{\hspace{1em}$\Box$\vspace{1em}}

\newcommand{\bC}{{\mathbb C}}
\newcommand{\bZ}{{\mathbb Z}}

\newcommand{\bQ}{{\mathbb Q}}
\newcommand{\bN}{{\mathbb N}}

\newcommand{\cA}{{\mathcal A}}

\newcommand{\rad}{{\frak r}}

\newcommand{\cE}{{\mathcal E}}

\newcommand{\kx}{{K[x]}}

\newcommand{\Ker}{\operatorname{\rm Ker\,}}
\newcommand{\im}{\operatorname{Im}}

\newcommand{\I}{\operatorname{I}}

\begin{document}

\title[Images of Ideals under Derivations and $\cE$-Derivations]
{Images of Ideals under Derivations and $\cE$-Derivations of  
Univariate Polynomial Algebras  over a Field  of Characteristic Zero}

  \author{Wenhua Zhao}      
    \date{\today}
\address{Department of Mathematics, Illinois State University, Normal, IL 61761. Email: wzhao@ilstu.edu}

\begin{abstract}
Let $K$ be a field of characteristic zero and $x$ a free variable. 
A $K$-$\cE$-derivation of $\kx$ is a $K$-linear map of the form $\I-\phi$ for 
some $K$-algebra endomorphism $\phi$ of $\kx$, 
where $\I$ denotes the identity map of $\kx$. 
In this paper we study the image of an ideal of $\kx$ under some  
$K$-derivations and $K$-$\cE$-derivations of $\kx$. 
We show that  
the LFED conjecture proposed in \cite{Open-LFNED} holds for all $K$-$\cE$-derivations 
and all locally finite $K$-derivations of $\kx$. We also show that  
the LNED conjecture proposed in \cite{Open-LFNED} holds for all locally 
nilpotent $K$-derivations of $\kx$, and also for all locally 
nilpotent $K$-$\cE$-derivations of $\kx$ and the ideals $u\kx$ such that 
either $u=0$, or $\deg u\le 1$, or $u$ has at least one repeated root 
in the algebraic closure of $K$.  
As a bi-product, the homogeneous Mathieu subspaces (Mathieu-Zhao spaces) 
of the univariate polynomial algebra 
over an arbitrary field have also been classified. 
\end{abstract}

   



%




%
%



\thanks{The author has been partially supported 
by the Simons Foundation grant 278638}

 \bibliographystyle{alpha}
    \maketitle


\renewcommand{\theequation}{\thesection.\arabic{equation}}
\renewcommand{\therema}{\thesection.\arabic{rema}}
\setcounter{equation}{0}
\setcounter{rema}{0}
\setcounter{section}{0}

\section{\bf Introduction}\label{S1}

Let $K$ be a field and $\cA$ a commutative $K$-algebra. 
We denote by $1_\cA$ or simply $1$ the identity element of $\cA$, if $\cA$ is unital, and 
$\I_\cA$ or simply $\I$ the identity map of $\cA$, if $\cA$ is clear in the context.  

A $K$-linear endomorphism $\eta$ of $\cA$ is said to be {\it locally nilpotent} (LN) 
if for each $a\in \cA$ there exists $m\ge 1$ such that $\eta^m(a)=0$, 
and {\it locally finite} (LF) if for each $a\in \cA$ the $K$-subspace spanned 
by $\eta^i(a)$ $(i\ge 0)$ is finite dimensional over $K$.    

A {\it $K$-derivation} $D$ of $\cA$ is a $K$-linear map 
$D:\cA \to \cA$ that satisfies $D(ab)=D(a)b+aD(b)$ for all $a, b\in \cA$. 
A {\it $K$-$\cE$-derivation} $\delta$ of $\cA$ 
is a $K$-linear map $\delta:\cA \to \cA$ such that for all 
$a, b\in \cA$ the following equation holds:
\begin{align}\label{ProdRule2}
\delta(ab)=\delta(a)b+a\delta(b)-\delta(a)\delta(b).  
\end{align}

It is easy to verify that $\delta$ is an $R$-$\cE$-derivation of $\cA$, 
if and only if $\delta=\I-\phi$ for some $R$-algebra endomorphism 
$\phi$ of $\cA$. Therefore an $R$-$\cE$-derivation is  
a special so-called $(s_1, s_2)$-derivation introduced 
by N. Jacobson \cite{J} and also a special 
semi-derivation introduced 
by J. Bergen in \cite{Bergen}. $R$-$\cE$-derivations have also been 
studied by many others under some different names such as 
$f$-derivations in \cite{E0, E} and 
$\phi$-derivations in \cite{BFF, BV}, etc..

 
Next, we recall the following two notions of associative algebras that were  
introduced in \cite{GIC, MS}. Since all algebras in this paper  
are commutative, here we recall only the cases for commutative algebras 
over a field.

\begin{defi} \label{Def-MS}
Let $K$ be a field and $\cA$ a commutative $K$-algebra. A $K$-subspace  
$V$ of $\cA$ is said to be a {\it Mathieu subspace} (MS)  
of $\cA$ if for all $a, b\in \cA$ with $a^m\in V$ for all $m\ge 1$, we have 
 $a^mb \in V$ for all $m\gg 0$. 
\end{defi}

Note that a MS is also called a {\it Mathieu-Zhao space} 
in the literature (e.g., see \cite{DEZ, EN, EH}, etc.), 
as suggested by A. van den Essen \cite{E2}. 
 
The introduction of this notion  
is mainly motivated by the study in \cite{Ma, IC} of 
the well-known Jacobian conjecture (see \cite{K, BCW, E}). 
See also \cite{DEZ}. But, a more interesting aspect 
of the notion is that it provides a natural but highly non-trivial    
generalization of the notion of ideals.

\begin{defi} \cite[p.\,247]{MS} \label{Def-Rad}
Let $V$ be a $K$-subspace (or a subset) of a $K$-algebra $\cA$. We define the {\it radical} 
$\rad(V)$ of $V$ to be 
\begin{align}
\rad(V)\!:=\{ a\in \cA\,|\, a^m\in V \text{ for all } m\gg0\}.
\end{align}
\end{defi}

%
%

Next we recall the cases of the so-called LFED and LNED conjectures 
proposed in \cite{Open-LFNED} for commutative algebras. 
For the study of some other cases of these two 
conjectures, see \cite{EWZ}, \cite{Open-LFNED}--\cite{LaurentPolyCase}.

\begin{conj}\label{LFED-Conj}
Let $K$ be a field of characteristic zero, $\cA$ a commutative $K$-algebra 
and $\delta$ a LF (locally finite) $K$-derivation or a LF 
$K$-$\cE$-derivation of $\cA$. Then the image 
$\im\delta\!:=\delta(\cA)$ of $\delta$ is a MS of $\cA$.   
\end{conj}

\begin{conj}\label{LNED-Conj}
Let $K$ be a field of characteristic zero, $\cA$ a commutative 
$K$-algebra and $\delta$ a LN (locally nilpotent) $K$-derivation or a LN 
$K$-$\cE$-derivation of $\cA$. Then $\delta$ maps every 
ideal of $\cA$ to a MS of $\cA$.  
\end{conj}

Throughout the paper we refer Conjecture \ref{LFED-Conj} as 
{\it the (commutative) LFED conjecture}, and Conjecture \ref{LNED-Conj} 
{\it the (commutative) LNED conjecture}.

In this paper, among some other results, we show   
the following two theorems regarding the commutative LFED and 
LNED conjectures, respectively.

\begin{theo}\label{MainThm-1}
Let $K$ be a field of characteristic zero and $x$ a free variable. 
Let $\delta$ be an arbitrary $K$-derivation or $K$-$\cE$-derivation of $\kx$.  
Then $\im\delta$ is a MS of $\kx$. 
In particular, the LFED conjecture \ref{LFED-Conj} 
holds for $\kx$. 
\end{theo} 
%
%

\begin{theo}\label{MainThm-2}
Let $K$ be a field of characteristic zero, $I$ an ideal of $\kx$ 
and $\delta$ a $K$-derivation or $K$-$\cE$-derivation of $\kx$. 
Then $\delta I$ is a MS of $\kx$ if one of the following conditions holds:
\begin{enumerate}
  \item[$1)$] $\delta$ is a locally nilpotent $K$-derivation of $\kx$; 
  \item[$2)$] $\delta=\I-\phi$ for some $K$-algebra 
endomorphism $\phi$ of $\kx$ such that $\deg \phi(x)\ge 2$;
  \item[$3)$] $\delta=\I-\phi$ for some $K$-algebra 
endomorphism $\phi$ of $\kx$ which maps $x$ to $x+c$ $(c\in K)$, 
and $I=u\kx$ such that either $u=0$, or $\deg u\le 1$,  
or $u$ has at least one repeated root in the algebraic closure of $K$. 
\end{enumerate}
\end{theo} 

From Theorem \ref{MainThm-2} above and also its proof it is easy to see that 
the LNED conjecture \ref{LNED-Conj} is established for $\kx$ 
except for the case that {\it $\delta=\I-\phi$ for some $K$-algebra 
endomorphism $\phi$ of $\kx$ that maps $x$ to $x+c$ $(0\ne c\in K)$, 
and $I=u \kx$ such that  $\deg u\ge 2$ and  
$u$ has no repeated root in the algebraic closure of $K$.  }

Theorems \ref{MainThm-1} and \ref{MainThm-2} are shown case by case. 
For some cases, e.g., the $K$-$\cE$-derivation case of Theorem \ref{MainThm-1},  
and Theorem \ref{MainThm-2}, $2)$, etc., certain stronger results are actually proved. 
Furthermore, as a bi-product of the proof of Theorem \ref{MainThm-2}, $2)$  
all homogeneous MSs of the univariate polynomial algebra over 
a field of arbitrary characteristic also are classified 
(see Proposition \ref{EquivH-MS}).\\



{\bf Arrangement.} In Section \ref{S2} we first show 
that the image of every $K$-derivation of $\kx$ is a MS 
(see Lemma \ref{LFD-Lma}), and then show that 
every LN $K$-derivation $D$, i.e., $D=a\frac{d}{dx}$ $(a\in K)$, 
maps each ideal of $\kx$ to a MS of $\kx$ (see Proposition \ref{Sub-MainThm1}). 
Consequently, the $K$-derivation cases of Theorems \ref{MainThm-1} 
and \ref{MainThm-2} are established. 
We also give an example, Example \ref{xddxExample}, 
to show that the LN condition in the LNED conjecture \ref{LNED-Conj} 
can not be replaced by the LF condition.

In Section \ref{S3} we let $K$ be a field of arbitrary characteristic and 
show the LFED conjecture \ref{LFED-Conj}  
for all the $K$-$\cE$-derivations $\delta$ that are not LN 
(the LN case of the conjecture also holds and follows from 
Lemma \ref{Lma-4.1} in Section \ref{S4}). 
In subsection \ref{S3.1} we consider 
the case $\delta=\I-\phi$ for some $K$-algebra endomorphism $\phi$ of $\kx$ that  
maps $x$ to $qx$ with $q\in K$. As a bi-product of the proof for this case 
we also obtain a classification of all  homogeneous MSs of $\kx$ 
(see Proposition \ref{EquivH-MS}). In subsection \ref{S3.2} 
we consider the case $\delta=\I-\phi$ such that $\phi$  
maps $x$ to $w(x)$ with $\deg w(x)\ge 2$. In particular, we show in 
Proposition \ref{ge2E-Case} that 
$\delta$ in this case actually maps every $K$-subspace 
to a MS of $\kx$, even though $\delta$ itself in this case is not LF (nor LN). 

In Section \ref{S4} we consider all LN $K$-$\cE$-derivations $\delta$ of $\kx$, i.e., $\delta=\I-\phi$ for some $K$-algebra endomorphism $\phi$ of $\kx$ that 
maps $x$ to $x+c$ $(c\in K)$, and prove Theorem \ref{MainThm-2}, $3)$.  Among all the cases studied in the paper, the proof of this case is the most involved, in which the Bernoulli polynomials; the Bernoulli numbers; and the 
Clausen-von Staudt Theorem \ref{Cl-St-Thm} that was found 
independently by Thomas Clausen \cite{Cl} and Karl von Staudt \cite{St} in  
$1840$, all unexpectedly play some crucial roles.

\section{\bf The Case of $K$-Derivations of $\kx$} \label{S2}

Throughout this section $K$ stands for a field of characteristic zero  
and $x$ a free variable. 
We denote by  $\p$  the $K$-derivation 
$d/d x$ of the univariate polynomial algebra $\kx$. 


Now, let $D$ be a nonzero $K$-derivation of $\kx$. Then 
$D=a(x)\p$ for some $a(x)\in \kx$. It is easy to 
see that $D$ is LF (locally finite), if and only if 
$\deg a\le 1$, and $D$ is LN (locally nilpotent), 
if and only if $\deg a=0$.  

We first show Theorem \ref{MainThm-1} for all 
$K$-derivations of $\kx$. 

\begin{lemma}\label{LFD-Lma}
Let $D=a(x)\p$ with $a(x)\in \kx$. Then 
$\im D=a(x)\kx$. 
In particular, the LFED conjecture \ref{LFED-Conj} 
holds for all $K$-derivations of $\kx$.
\end{lemma}

\pf Since $\im \p=K[x]$, 
we have $\im D=a(x)\kx$, which is  
an ideal of $\kx$, and hence 
also a MS of $\kx$.
\epfv

%

Next, we consider LN (locally nilpotent) $K$-derivations. First, 
let us recall the following:  

\begin{theo}\label{Int-ab-Thm}
Let $a \ne b\in K$ and set 
\begin{align} \label{V-ab}
V_{a, b}\!:=\left\{f\in  \kx\,\big|\, \int_a^b f(x) dx=0\right\}. 
\end{align}
Then  $\rad(V_{a, b})=0$.
\end{theo} 

For an algebraic proof of the theorem above, see \cite[Theorem $4.1$]{FPYZ}, 
and for a complex analytic proof with some slightly stronger condition, 
see \cite[Corollary $4.3$]{P}. Although the theorem is proven in \cite{FPYZ} and \cite{P} over the complex field $\bC$, by the Lefschetz rule or from the proof of \cite[Theorem $4.1$]{FPYZ} it is easy to see that the theorem also holds over all the fields 
of characteristic zero.

%

Now we show Theorem \ref{MainThm-2} for LN  
$K$-derivations $D$ of $\kx$. Since the only LN 
$K$-derivations of $\kx$ are $a\p$ with $a\in K$, 
it suffices to consider the case that $D=\p$, 
for the case $D=0$ is trivial.    
 
\begin{propo}\label{Sub-MainThm1}
Let $I$ be a nonzero ideal of $\kx$. Write $I=(u(x))$ 
for some $u(x)\in \kx$. Then the following statements hold:   
\begin{enumerate}
  \item[$1)$] if $u(x)=(x-c)^n$ for some $c\in K$ and $n\ge 0$, 
  then the image $\p I$ of $I$ under $\p$ 
  is an ideal of $\kx$. More precisely, 
\begin{align}
\p I=\begin{cases}
\kx &\text{ if } n\le 1;\\
 (x-c)^{n-1}\kx &\text{ if } n\ge 2.
\end{cases}
\end{align}
\item[$2)$] if $u(x)\ne (x-c)^n$ for any $c\in K$ and $n\ge 0$, then 
the radical $\rad(\p  I)=\{0\}$, whence $\p I$ is a MS of $\kx$. 
\end{enumerate}
Consequently, the LNED conjecture \ref{LNED-Conj} 
holds for all LN $K$-derivations of $\kx$.
\end{propo}
 
\pf Since $(x-c)^n$ $(n\ge 0)$ form a $K$-linear basis of $\kx$, 
it is easy to see that statement $1)$ holds.

To show statement $2)$, note first that by the assumption we have that  
$\deg u\ge 2$ and has at least two distinct roots $a$ and $b$ 
in the algebraic closure $\bar K$ of $K$. Now for every $f\in \p I$,  
write $f=\p(ug)$ for some $g\in \kx$. Then we have 
$\int_a^b f dx=\int_a^b \p(ug) dx=(ug)|_a^b=0$.
Therefore $\p I \subseteq V_{a, b}\!:=\{f\in \bar K[x]\,|\, \int_a^b fdx=0\}$.  
Applying Theorem \ref{Int-ab-Thm} to $\bar K[x]$ we get $\rad( V_{a, b})=\{0\}$, 
whence $\rad(\p I)=\{0\}$ and statement $2)$ follows. 
\epfv
 
It is worthy to point out that the LNED conjecture \ref{LNED-Conj} 
can not be generalized to all LF $K$-derivations of $\kx$, 
which can be seen from the following:  
 
\begin{exam}\label{xddxExample}
Let $D=x\p$ and $I=(x^2-1)\kx$.  
Then the image $D I$ of $I$ under $D$ is not a MS 
of $\kx$. 
\end{exam}

\pf Let $V=D I$. Then for all $k\ge 0$, since  
$D(x^{k+2}-x^k)=(k+2)x^{k+2}-kx^k \in V$,  
we have 
\begin{align}\label{JumpbyTwo}
(k+2)x^{k+2}\equiv kx^k \mod V. 
\end{align}
In particular, $x^2\in V$ (by letting $k=0$) and, inductively by 
Eq.\,(\ref{JumpbyTwo}), so are $x^{2n}$ for all $n\ge 1$.

On the other hand, $x\not \in V$ since each nonzero element of $V$ has degree 
at least $2$, and by Eq.\,(\ref{JumpbyTwo}) neither are $x^{2n+1}$ 
for all $n\ge 1$. Therefore $(x^2)^mx\not \in V$ for all $m\ge 1$, 
whence $V$ is not a MS of $\kx$.  
\epfv
 
We end this section with the following remark on an application 
of the results proved in this section.

\begin{rmk}
Let $u, v \in \kx$, $D=u\p$,  $I=v(x)\kx$ and $\Lambda$ 
the differential operator of $\kx$ that maps $f\in\kx$ 
to $u\p(vf)$, i.e., 
$$
\Lambda\!:=u(v\p+v')
$$
Then it is easy to see that $\im \Lambda=DI$. 
Therefore by Lemma \ref{Lma-4.1} and Proposition \ref{Sub-MainThm1} 
we see that $\im \Lambda$ is a MS of $\cA$, if $\deg u\le 1$. 
For example, by letting $u=1$ we see that for all  $v(x)\in \kx$ 
the image of the differential operators $v(x)\p+v'(x)$ is a MS 
of $\kx$.  

For some other differential operators with the image 
being a MS, see \cite{IC, GIC, EZ, EWZ}.  
\end{rmk}

\section{\bf The Case of $\cE$-Derivations of $\kx$}\label{S3}

Through this section {\bf $\mathbf{K}$ denotes a field of arbitrary characteristic}, and 
$\phi$ a $K$-algebra endomorphism of $\kx$, and 
$\delta=\I-\phi$. Since the case $\phi=0$ is trivial, we assume $\phi\ne 0$. 
Denote by $w(x)$ the image of $x\in \kx$ 
under $\phi$, i.e., $w(x)=\phi(x)$. 
Then $\phi$ is completely determined 
by $w(x)$. More precisely, for each $f(x)\in \kx$, 
we have $\phi(f)=f(w(x))$.  
 
We start with the following 

\begin{lemma}\label{Lemma3.1}
Assume $w(x)=ax+b$ with $a, b\in K$ and $a\ne 0, 1$. 
Let $\psi$ be the $K$-algebra automorphism of $\kx$ which maps $x$ to $x+(1-a)^{-1} b$. Then 
$\psi\circ\phi\circ\psi^{-1}$ is the $K$-algebra automorphism 
of $\kx$ which maps $x$ to $ax$.   
\end{lemma}

\pf Note that the inverse map $\psi^{-1}$ is the $K$-algebra automorphism of $\kx$ which maps $x$ to $x-(1-a)^{-1}b$. Then   we have  
\begin{align*}
\psi \circ \phi\circ \psi^{-1}  (x) 
&=(\psi\circ \phi) ( x-(1-a)^{-1}b\big )=
\psi  (ax+b-(1-a)^{-1}b)\\
&= a(x+(1-a)^{-1}b)+b-(1-a)^{-1}b \\
&= ax+b+(a(1-a)^{-1}-(1-a)^{-1})b \\
&=ax+b-b=ax.
\end{align*}
Hence the lemma follows.
\epfv

Note that $K$-algebra automorphisms preserve ideals and MSs, and conjugations by 
$K$-algebra automorphisms preserve (LF or LN) derivations and $\cE$-derivations. 
By the lemma above the proofs of Theorem \ref{MainThm-1} and \ref{MainThm-2} 
for $K$-$\cE$-derivations of $\kx$ can be divided into the following 
four (exhausting) cases:

{\it
\begin{enumerate}\label{CaseList}
  \item[I)] $\deg w=0$, i.e., $w(x)=c$ for some $c\in K$;
  \item[II)] $w(x)=x+c$ for some $c\in K$;
  \item[III)] $w(x)=qx$ for some nonzero $q\in K$;
  \item[IV)] $\deg w(x)\ge 2$.
  \end{enumerate}
}

For Case I it is easy to verify, or by the more 
general \cite[Proposition 5.2]{Open-LFNED}, 
that we have the following: 

\begin{lemma}\label{CaseI-Lma}
Let $c\in K$ and $\phi$ the $K$-algebra endomorphism that maps 
$f\in\kx$ to $f(c)$. Then the image $\im (\I-\phi)=\Ker \phi=(x-c)\kx$,  
and hence is a MS of $\kx$. 
\end{lemma}
 
%
%
%
%

%

In the rest of this section we consider Case III in subsection \ref{S3.1} 
and Case IV in subsection \ref{S3.2}.  \\

\subsection{{\bf Case III}: $\pmb{w(x)=qx}$} \label{S3.1}
Note that this case has been shown in \cite[Corollary 3.15]{LaurentPolyCase} 
for multivariate polynomial algebras over a field of characteristic zero.
Here we give a more straightforward proof over the field $K$ 
(of arbitrary characteristic). As a bi-product the homogeneous MSs 
of $\kx$ are also classified.

\begin{lemma}\label{qx-Case}
Let $0\ne q\in K$ and $\phi$ the $K$-algebra endomorphism of $\kx$ that maps $x$ to $qx$. 
Set $\delta=\I-\phi$. Then the following statements hold: 
\begin{enumerate}
\item[$1)$] If $q=1$, then $\im\delta=0$.

\item[$2)$] If $q$ is not a root of unity in $K$, 
then $\im\delta=x\kx$. 

\item[$3)$]  If $q$ is a root of unity in $K$, 
then $\rad(\im\delta)=\{0\}$.
\end{enumerate}
In all the cases above $\delta I$ is a MS of $\kx$.
\end{lemma}

\pf $1)$ In this case $\phi=\I$ and $\delta=0$. Hence $\im\delta=0$.

$2)$ For all $n\ge 1$,   
we have 
\begin{align}\label{qx-Case-peq1}
\delta x^n=x^n-(qx)^n=(1-q^n)x^n.
\end{align} 
Since $q$ is not a root of unity, we have $1-q^n\ne 0$, and hence $x^n\in \im \delta$,  
for all $n\ge 1$. Since $\delta 1=0$, we have $\im \delta =x\kx$, i.e., 
statement $2)$ follows.   

$3)$ If $q=1$, the statement follows from statement $1)$. Assume $q\ne 1$ and 
let $r$ be the order of $q$, i.e., the least positive  
integer such that $q^r=1$. Then $r\ge 2$. 
Let $\{n_i\,|\,i\ge 1\}$ be the sequence of 
all positive integers $n$ such that $r \nmid n$. Then 
by Eq.\,(\ref{qx-Case-peq1}) it is easy to see that 
$\im \delta$ is the homogeneous $K$-subspace spanned   
by the monomials $x^{n_i}$ $(i\ge 1)$ over $K$. Note that for each integer   
$d\ge 1$, we have $dr \not \in \{n_i\,|\,i\ge 1\}$. 
Then the statement immediately follows from the lemma below.
\epfv
 
\begin{lemma}\label{EquivH-MS}
Let  
$\{ n_i\,|\, i\ge 1 \}$ be a strictly increasing (infinite)  
sequence of positive integers such that $n_{i+1}-n_i\ne 1$ 
for infinitely many $i\ge 1$. Let $V$ be the (homogeneous) $K$-subspace of 
$K[x]$ spanned by $x^{n_i}$ $(i\ge 1)$ over $K$. 
Then the following three statements are equivalent:
\begin{enumerate}
  \item[$1)$] $\rad(V)=\{0\}$;
  \item[$2)$] $V$ is a MS of $\kx$;
  \item[$3)$] there exists no integer $d\ge 1$ such that 
 $md\in \{n_i\,|\, i\ge 1\}$ for all $m\ge 1$. 
\end{enumerate}    
\end{lemma}
\pf $1)\Rightarrow 2)$ is obvious. 

$2)\Rightarrow 3)$: Assume otherwise. Let $d\ge 1$ such that  
$md\in \{n_i\,|\, i\ge 1\}$ for all $m\ge 1$. Hence  
$x^{md}\in V$ for all $m\ge 1$. 
If $d=1$, then the sequence 
$\{n_i\,|\, i\ge 1\}$ contains all positive integers, which 
contradicts to the assumption on the sequence 
$\{n_i\,|\, i\ge 1\}$. 
So we have $d\ge 2$. 

Since $(x^d)^m=x^{dm}\in V$ for all $m\ge 1$ and $V$ is a MS of $\kx$,  
for each $0\le r\le d-1$ 
there exists $N_r\ge 1$ such that for all $m\ge N_r$, we have $x^{md+r}=(x^d)^m x^r \in V$, 
and hence $md+r\in \{n_i\,|\, i\ge 1\}$. 
Let $N=\max\{N_r\,|\,0\le r\le d-1\}$. Then for all   
$k\ge Nd$ we have $x^k\in V$, whence 
$k\in \{n_i\,|\,i\ge 1\}$, which contradicts 
again to the assumption on the sequence 
$\{n_i\,|\, i\ge 1\}$.  
Hence statement $3)$ follows.

$3)\Rightarrow 1)$: Assume otherwise. Let 
$0\ne f(x)\in \rad(V)$, i.e., 
$f^m(x)\in V$ when $m\gg 0$.  Since $1\not \in V$, 
we have $d\!:=\deg f(x)\ge 1$, and since $V$ is homogeneous, 
we further have $x^{dm}\in V$. Hence $md\in \{n_i\,|\, i\ge 1\}$, 
for all $m\gg 0$. Replacing $d$ by a multiple of $d$ we have 
$md\in \{n_i\,|\, i\ge 1\}$ for all $m\ge 1$, which contradicts 
to statement $3)$. 
\epfv

One bi-product of the lemma above is 
the following classification of 
homogeneous MSs of univariate polynomial 
algebra over an arbitrary field. 

\begin{propo}\label{cls-H-MS}
Let  
$V$ be a homogeneous subspace of $\kx$. Then $V$ 
is a MS of $\kx$, if and only if one of the 
following conditions holds:
\begin{enumerate}
\item[$1)$]  $V=\kx$; 
\item[$2)$]  $\dim_K V<\infty$ and $1\not\in V$;  
\item[$3)$]  $1\not\in V$ and there exists $N\ge 1$ such that $(x^N) \subseteq V$;   
\item[$4)$]  $V$ is spanned by $x^{n_i}$ over $K$ 
for a strictly increasing (infinite)  
sequence of positive integers $\{n_i\,|\, i\ge 1\}$  such that 
there exists no integer $d\ge 1$ with $md\in \{n_i\,|\, k\ge 1\}$ for all $m\ge 1$. 
\end{enumerate} 
\end{propo}

\pf $(\Leftarrow)$ First, if $V$ satisfies statements $1)$, $2)$ or $3)$, then 
it is easy to check directly by Definition \ref{Def-MS} that $V$ is indeed a MS of $\cA$. 
If $V$ satisfies statement $4)$, then it is easy to see that $n_{i+1}-n_i\ne 1$ 
for infinitely many $i\ge 1$, and by Lemma \ref{EquivH-MS} $V$ is a MS of $\kx$. 

$(\Rightarrow)$ Assume that statement $1)$ fails, i.e., $V\ne \kx$. 
Consider the case $\dim_K V<\infty$. If $1\in V$, then 
by Definition \ref{Def-MS} we have $V=\kx$, Contradiction. 
Hence in this case statement $2)$ holds. 

Consider the case $\dim_K V=\infty$. If $1\in V$, then by Definition \ref{Def-MS} 
we have $V=\kx$, contradiction again. So $1\not \in V$ and  
there exists an infinite increasing sequence 
$\{n_i\,|\, i\ge 1\}$ of positive integers such that 
$V$ is spanned over $K$ by $x^{n_i}$ $(i\ge 1)$. 

If statement $3)$ does not hold, 
then there are infinitely many $i\ge 1$ 
such that $n_{i+1}-n_i\ne 1$, and 
by Lemma \ref{EquivH-MS} statement $4)$ holds.
\epfv

%
%
%

Next, we give the following example to show that 
the LNED conjecture \ref{LNED-Conj}  can not be generalized 
to all LF $\cE$-derivations.

\begin{exam}\label{qxExample}
Let $0\ne q\in K$, $\phi$ the $K$-algebra endomorphism of $\kx$ that maps $x$ to $qx$,  
and $I$ the ideal of $\kx$ generated by $x^2-1$. Assume that $q$ is not a root of 
unity. Set $\delta\!:=\I-\phi$. Then $\delta$ is LF but the image $\delta I$ of $I$ 
under $\delta$ is not a MS of $\kx$. 
\end{exam}

\pf Note first that for all $k\ge 0$, we have
$$
\delta(x^{k+2}-x^k)=(1-q^{k+2})x^{k+2}-(1-q^k)x^k \in \delta I.
$$ 
Hence  
\begin{align}\label{q-JumpbyTwo}
(1-q^{k+2})x^{k+2}\equiv (1-q^k)x^k  \mod \delta I. 
\end{align}
In particular, by letting $k=0$ we have $x^2\in \delta I$, for $q$ is not root of unity.  
Then by Eq.\,(\ref{q-JumpbyTwo}) inductively $x^{2n}\in \delta I$ for all $n\ge 1$. 

On the other hand, $x\not \in \delta I$, for $q$ is not root of unity and 
each nonzero element of $\delta I$ has the degree at least $2$.  
Then by Eq.\,(\ref{q-JumpbyTwo}) neither are $(x^2)^mx=x^{2m+1}$ 
for all $m\ge 1$. Hence $\delta I$ is not a MS of $\kx$.  
\epfv

\subsection{{\bf Case  IV:}  $\pmb{\deg w(x)\ge 2}$} \label{S3.2}
In this subsection we fix a $K$-algebra endomorphism $\phi$ of $\kx$ that 
maps $x$ to $w(x)$ with $d\:=\deg w(x)\ge 2$. Set $\delta\!:=\I-\phi$. 
Write $w(x)=\sum_{i=0}^d a_i x^i$ with $d\ge 2$, $a_d\ne 0$ and all $a_i$'s in $K$.  
Then Theorem \ref{MainThm-2}, 2) for $\delta$ immediately follows from the following: 

\begin{propo}\label{ge2E-Case}
Let $f\in \kx$ such that $f^i\in \im\delta$ for all $1\le i\le 3$. 
Then $f=0$.
Consequently, $\rad(\im\delta)=0$ and $\delta$ maps 
every $K$-subspace of $\kx$ to  
a MS of $\kx$.
\end{propo}

To prove the proposition above we first fixed the 
following notations.

Let $W=K[w(x)]$,  
$\Lambda=\bN\backslash d\bN$ and 
$U$ be $K$-subspace of $\kx$ spanned by $x^{m}$ 
($m\in \Lambda$) over $K$. 
Since $d=\deg w\ge 2$, it is easy to see that 
for each $f\in \kx$ there exist  
unique $f_1\in U$ and $f_2\in W$ such that 
$f=f_1+f_2$. In this case we set 
$\ell(f)=\deg f_1$ if $f_1\ne 0$, and $\ell(f)=0$, 
otherwise. 

With the setting above the following lemma (with assumption 
$d\ge 2$) can be easily verified. 

\begin{lemma}\label{Lma-2.2.1}
$1)$ $\deg f \ge \ell(f)$ for all nonzero $f\in\kx$. 

$2)$ For all $f, g\in \kx$ with $f\equiv g\mod W$ we have 
$\ell(f)=\ell(g)$.
\end{lemma}
 
We also need the following two lemmas.  

\begin{lemma}\label{Lma-2.2.1-2}
Let $0\ne f\in \kx$ with $d \nmid  \deg f$. Then 
$\ell(f)=\deg f$.  
\end{lemma}

\pf Assume otherwise. Then there exists $g\in U$ such that 
$f-g\in W$ and $\deg f\ne \ell (f)=\deg g$. Then $f-g\ne 0$ and 
$\deg (f-g)$ is a multiple of $d$. On the other hand, $\deg (f-g)$
is equal to either $\deg f$ or $\deg g$. Therefore, 
either $\deg f$ or $\deg g$ is a multiple of $d$, 
which is a contradiction.
\epfv

\begin{lemma}\label{Lma-2.2.2}
For all nonzero $f\in \im\delta$, the following statements hold:
\begin{enumerate}
  \item[$1)$] if $f\in W$, say $f(x)=\tilde f(w(x))$, then $\tilde f\in \im \delta$;
  \item[$2)$] if $f\not\in W$, then $\deg f \ge d\ell(f)\ge d$. 
\end{enumerate}
\end{lemma}

\pf $1)$ Write $f(x)=\delta u=u(x)-u(w(x))$ for some $u\in \kx$. 
Then  
$\tilde f(w)=u(x)-u(w)$ and $u(x)=\tilde f(w)+u(w)$. Hence 
 $u(x)=\tilde u(w)$ with $\tilde u=\tilde f+u$, and  
$$
\tilde f(x)=\tilde u(x)-u(x)=\tilde u(x)-\tilde u(w)=\delta \tilde u.
$$
 Therefore $\tilde f\in \im\delta$.

$2)$ Since $f\not \in W$, we have $f\ne 0$ and $\ell(f)\ge 1$.  
Write $f(x)=u(x)-u(w(x))$ for some nonzero $u\in \kx$. Since $d=\deg w\ge 2$, 
we have $\deg f=d \deg u$. By Lemma \ref{Lma-2.2.1} we also have 
$\deg u\ge \ell(u)$ and $\ell(f)=\ell(u)$. Therefore  we have 
$$
\deg f= d\deg u \ge d\ell(u)=d\ell(f)\ge d. 
$$
\epfv

\underline{\bf Proof of Proposition \ref{ge2E-Case}}: Assume otherwise, i.e., 
$f\ne 0$. If $f\in W$, say $f=\tilde f(w)$, then $f^2, f^3\in W$ and by 
applying lemma \ref{Lma-2.2.2}, $1)$ to $f^i$ $(1\le i\le 3)$  
we have  $\tilde f^i\in \im\delta$ for all $1\le i\le 3$. 
Since $\deg f=d\deg \tilde f>\deg \tilde f$, by repeating the procedure that replaces 
$f$ by $\tilde f$, whenever it is possible, we may assume that 
$f\not\in W$ and $f^i \in \im\delta$ $(1\le i\le 3)$. In particular, $f\ne 0$.
Furthermore, since $d=\deg w\ge 2$, $\im \delta$ obviously does not contain any nonzero 
constant polynomials. Therefore we also have $\deg f \ge 1$.

Write $f(x)=u(x)-u(w)$ for some nonzero $u(x)\not \in W$ 
with $\deg u(x)\ge 1$. 
Then $\ell(u)\ge 1$ and 
\begin{align}\label{ge2E-Case-peq0}
\deg f=d\deg u.
\end{align}

Assume first that char.\,$K=0$ or $p>2$. Then 
applying Lemma \ref{Lma-2.2.1}, $1)$ to 
$f^2$ we have 
\begin{align}\label{ge2E-Case-peq1}
2\deg f\ge d\ell(f^2).
\end{align}

On the other hand, we also have 
$$
f^2=(u(x)-u(w))^2=u^2(x)-2u(x)u(w)+u^2(w).
$$ 
Then by Lemma \ref{Lma-2.2.1}, 2) and Lemma \ref{Lma-2.2.1-2} as well as 
the assumption $d=\deg w\ge 2$ we get  
\begin{align}\label{ge2E-Case-peq3}
\ell(f^2)=\deg u(x)u(w)=(d+1)\deg u.
\end{align}
Combining Eqs.\,(\ref{ge2E-Case-peq0})--(\ref{ge2E-Case-peq3}) 
we have 
$$
2\deg f \ge d(d+1)\deg u=(d+1)\deg f, 
$$
which is a contradiction, for $d\ge 2$ and $\deg f\ge 1$.   

Now assume char.\,$K=2$. Then 
applying Lemma \ref{Lma-2.2.1}, $1)$ to 
$f^3$ we have 
\begin{align}\label{ge2E-Case-peq4}
3\deg f\ge d\ell(f^3).
\end{align}

On the other hand, we also have 
$$
f^3=(u(x)-u(w))^3=u^3(x)+u^2(x)u(w)+u(x)u^2(w)+u^3(w).
$$ 
Then by Lemma \ref{Lma-2.2.1}, 2) and Lemma \ref{Lma-2.2.1-2} as well as 
the assumption $d=\deg w\ge 2$ we get  
\begin{align}\label{ge2E-Case-peq5}
\ell(f^3)=\deg u(x)u^2(w)=(2d+1)\deg u.
\end{align}
Combining Eqs.\,(\ref{ge2E-Case-peq0}), (\ref{ge2E-Case-peq4}) and 
(\ref{ge2E-Case-peq5}) 
we have 
$$
3\deg f \ge d(2d+1)\deg u=(2d+1)\deg f, 
$$
which is a contradiction, for $d\ge 2$ and $\deg f\ge 1$.  
Hence the proposition follows.   
\epfv

\begin{rmk}
By Lemmas \ref{CaseI-Lma} and \ref{qx-Case}, Proposition \ref{ge2E-Case} and also 
Lemma \ref{Lma-4.1} in the next section we see that Theorem \ref{MainThm-1} holds for 
all $K$-$\cE$-derivations (not necessarily LF) of $\kx$. Furthermore, 
Proposition \ref{ge2E-Case} also implies statement $2)$ of 
Theorem \ref{MainThm-2}.
\end{rmk}

\section{\bf The Locally Nilpotent $\cE$-Derivation Case} 
\label{S4}
In this section we let $K$ be a field of characteristic zero. We  
consider the LN (locally nilpotent) $K$-$\cE$-derivations of $\kx$ and 
give a proof for statement $3)$ of Theorem \ref{MainThm-2}.  
From the exhausting list on page 
\pageref{CaseList} it is easy to see that 
the only nonzero LN $K$-$\cE$-derivations 
$\delta$ of $\kx$ are those in Case II of the list, 
i.e., $\delta=\I-\phi$, where $\phi$ is the affine 
translation that maps $x$ to $x+c$ for some $c\in K$.

Note that, if $c=0$, we have $\phi=\I$ and 
$\delta$ is the zero map, which is a trivial case. 
So throughout this section we assume $\delta=\I-\phi$, and 
$\phi$ maps $x$ to $x+c$ for a fixed $0\ne c\in K$. 

  
We first consider the images of the 
following ideals of $\kx$ under $\delta$. 

\begin{lemma}\label{Lma-4.1}
Let $\phi$, $\delta$ be as above and $I=\kx$ or $(x-a)\kx$  
for some $a\in K$. Then $\delta I =\kx$.  
\end{lemma}

\pf It suffices to show the case that $I=(x-a)$, for 
$\delta I\subseteq  \delta (\kx)$.

First, since $\delta (x-a)=-c$ and $c\ne 0$, we see that 
$1=x^0\subset \delta I$. Assume that for some $n\ge 1$ 
all polynomial $f(x)\in \kx$ with $\deg f\le n-1$ 
lie in $\delta I$. Consider $\delta\big((x-a)^{n+1}\big)\in \delta I$:    
$$
\delta\big((x-a)^{n+1}\big)=(x-a)^{n+1}-(x+c-a)^{n+1}
=-(n+1)cx^n+h(x) 
$$
for some $h(x)\in \kx$ with $\deg h\le n-1$. 
By the induction assumption above we have $h\in \delta I$. 
Therefore $-(n+1)cx^n\in \delta I$ and hence so does $x^n$. 
Therefore 
all polynomial $f(x)\in \kx$ with $\deg f\le n$ 
lie in $\delta I$, whence by induction the lemma follows. 
\epfv

In order to consider the images under $\delta$ of the 
ideals of $\kx$ generated by polynomials of the degree $\ge 2$,    
we first need to recall some well-known facts on the Bernoulli polynomials 
and the Bernoulli numbers (e.g., see \cite{Wiki}, \cite{Ber} 
and the references therein). 
   
First, {\it the Bernoulli polynomials} $\{B_n(t)\,| \, n\ge 0\}$ 
are defined by the following generating function:
\begin{align}
\frac{ue^{tu}}{e^u-1}=\sum_{n=0}^\infty B_n(t)\frac{u^n}{n!}.
\end{align}

For example, the first four Bernoulli polynomials are given 
as follows:  
\begin{align}
B_0(t)=1,\nno  
\end{align}
\begin{align}
B_1(t)=t-\frac12,\label{B1}  
\end{align}
\begin{align}
B_2(t)=t^2-t+\frac16,\nno 
\end{align}
\begin{align}
B_3(t)=t^3-\frac32 t^2+\frac12 t.\nno 
\end{align}

For every $n\ge 0$, the following 
identities of the Bernoulli polynomials hold:
\begin{align}
B_n(t+1)-B_n(t)=nt^{n-1};\label{BerId-0} 
\end{align}
\begin{align}
\frac{d}{dt}B_{n+1}(t)=(n+1)B_n(t).\label{BerId-2}  
\end{align}
\begin{align}
\sum_{k=0}^n\binom{n+1}{k}B_k(t)=(n+1)t^n, \label{BerId-1} 
\end{align}
\begin{align}
 (-1)^n B_n(-t) = B_n(t) + nt^{n-1}.  \label{BerId-5}
\end{align}

One remark on the Bernoulli polynomials is the following: 

\begin{propo}\label{BernoulliSubspace}
$1)$ the $K$-subspace spanned by the Bernoulli polynomials $B_n(x)$ $(n\ge 1)$ coincides with the $K$-subspace $V_{0, 1}$ defined in Eq.\,(\ref{V-ab}) with $a=0$ and $b=1$.

$2)$ Let $\Lambda$ be a non-empty set of positive integers and $W$ the 
$K$-subspace of $\kx$ spanned by the Bernoulli polynomials $B_i(x)$ with $i\in \Lambda$.  
Then $W$ is a MS of $\kx$ with $\rad(W)=\{0\}$.  
\end{propo}

\pf $1)$ By Eq.\,(\ref{BerId-0}) with $t=0$ we have $B_n(1)=B_n(0)$ for all $n\ge 2$. 
Then  by Eq.\,(\ref{BerId-2}) we have for all $n\ge 1$   
\begin{align}\label{BernoulliSubspace-peq1}
\int_0^1 B_n(x) dx=\frac{1}{n+1}B_{n+1}(x)\big|_0^1=0.
\end{align}
Hence $B_n(x)\in V_{0,1}$ for all $n\ge 1$.

Conversely, by Eqs.\,(\ref{BerId-2}) and the fact $B_0(x)=1$ we have 
$\deg B_n(x)=n$ for all $n\ge 0$. Hence the Bernoulli polynomials $B_n(x)$ $(n\ge 0)$
form a $K$-linear basis of $\kx$. In particular, every  
$f(x)\in\kx$ can be written uniquely as $f(x)=\sum_{i=0}^d c_iB_i(x)$ with 
$c_i$'s in $K$. Then by Eqs.\,(\ref{BernoulliSubspace-peq1}) and the fact $B_0(x)=1$ 
we see that $f(x)\in V_{0, 1}$, if and only if $c_0=0$. 
Hence statement $1)$ follows.
  
$2)$ By statement $1)$ we have $W\subseteq V_{0, 1}$. Then 
by Theorem \ref{Int-ab-Thm} we have 
 $\rad(W)=\{0\}$, whence $W$ is a MS of $\kx$. 
\epfv

The constant term $B_n$ of the Bernoulli polynomial $B_n(t)$ is called {\it the $n^\text{th}$  
Bernoulli number}, i.e., $B_n\!:=B_n(0)$. The Bernoulli polynomials $B_n(t)$ 
can be expressed in terms of the Bernoulli numbers $B_n$ $(n\ge 0)$ as follows:
\begin{align}
B_n(t)= \sum_{k=0}^n \binom{n}{k}B_k \,t^{n-k}.\label{BerId-4}
\end{align}



For the Bernoulli numbers we have the following remarkable theorem, which was found independently by   
Thomas Clausen (\cite{Cl}, 1840) and  Karl von Staudt (\cite{St}, 1840). 
See also \cite{Wiki2}.

\begin{theo}[{\bf The Clausen-von Staudt Theorem}]\label{Cl-St-Thm}
For each $n\ge 1$, we have
$$B_{2n}+\sum_{\substack{q:\,prime\\ (q-1)|2n}} \frac1q\in \bZ.$$
\end{theo}
%
%

From the Clausen-von Staudt Theorem we immediately have the following 

\begin{corol} \label{pv4Ber}
Let $p$ be an odd prime and $\nu_p(\cdot)$ the $p$-valuation on $\bQ$.
Then the following statements hold:  
\begin{enumerate}
\item[$1)$] $\nu_p (B_n)\ge 0$ for all $1\le n\le p-2$;
\item[$2)$] $\nu_p (B_{p-1})<0$, i.e., 
$p$ divides the denominator but not the numerator of the reduced fraction 
form of $B_{p-1}$. 
\end{enumerate}
\end{corol}

Now we get back to the image of an ideal $I$ under 
the fixed $K$-$\cE$-derivation $\delta$ of $\kx$. 

\begin{lemma}\label{Dn-Lma}
Let $I$ be the ideal of $\kx$ generated by the polynomial $x^2-a x$ for some 
$a\in K$ and $\beta\!:=a/c$. Set for all $n\ge 0$  
\begin{align}\label{Def-Dn}
D_n(t)\!:=\frac{B_{n+1}(t)-B_{n+1} }{(n+1)t}.  
\end{align}   
Then for all $n\ge 0$ we have $D_n(t)\in \bQ[t]$ and 
\begin{align}\label{Dn-Lma-eq2}
x^n\equiv D_n(\beta)c^n \mod \delta I.
\end{align}
\end{lemma}

Note that by Eqs.\,(\ref{BerId-0}) and (\ref{BerId-4}), $D_n(t)$ $(n\ge 0)$ are actually given by
\begin{align}
D_n(t)=\frac1{n+1}\sum_{i=0}^{n} \binom{n+1}{i}B_{i}\,t^{n-i}. \label{Dk-Lma-eq4}
\end{align}

\vspace{6mm}

\underline{\it Proof of Lemma \ref{Dn-Lma}}: First, for each $n\ge 1$, we have $x^{n+1}-a x^n\in I$, and hence 
$\delta x^{n+1}\equiv a \delta x^n  \mod \delta I$. Consequently, 
$\delta x^{n+1} \equiv a^n \delta x \mod \delta I$. 
Since $\delta x=-c$,  we have
$$
(x+c)^{n+1}-x^{n+1} \equiv a^n c \mod \delta I.
$$
$$
\sum_{k=0}^n \binom{n+1}k c^{n+1-k}x^k \equiv a^n c \mod \delta I. 
$$ 
\begin{align}\label{Lma-3.9-peq-3}
\sum_{k=0}^n \binom{n+1}k c^{-k} x^k \equiv a^n/c^n =\beta^n \mod \delta I. 
\end{align}

From the equation above it is easy to see recursively that 
for each $k\ge 0$, we have $x^k \equiv E_k c^{k}$ for some $E_k\in K$, 
which is a polynomial in $c$, $c^{-1}$ and $a$ with coefficients in $\bQ$. 
Furthermore, since $\delta I$ obviously does not contain any nonzero 
constant, $E_k$ $(k\ge 0)$ are actually unique. In particular, $E_0=1$.


Now we plug the relations $x^k \equiv   E_k  c^{k} \mod \delta I$ $(k\ge 1)$ 
into Eq.\,(\ref{Lma-3.9-peq-3}) and get
\begin{align}\label{Lma-3.9-peq-4}
\sum_{k=0}^n \binom{n+1}k E_k \equiv \beta^n \mod \delta I.
\end{align} 

Since $\delta I$ does not contain any nonzero constant, the equation above 
is the same as 
\begin{align}\label{Lma-3.9-peq-5}
\sum_{k=0}^n \binom{n+1}k E_k = \beta^n.
\end{align}

On the other hand, taking $\int_0^t$ to Eq.\,(\ref{BerId-1}) and applying 
Eq.\,(\ref{BerId-2}) we get
\begin{align}\label{Lma-3.9-peq-6}
\sum_{k=0}^n\binom{n+1}{k}\frac{B_{k+1}(t)-B_{k+1}}{k+1}=t^{n+1}.   
\end{align}

Since $B_{k+1}=B_{k+1}(0)$, $B_{k+1}(t)-B_{k+1}$ 
is divisible by $t$.  Then by Eq.\,(\ref{Def-Dn}) we have $D_k(t)\in \bQ[t]$, 
for $B_k(t)\in\bQ[t]$, for all $k\ge 0$. Furthermore, 
the equation above can be re-written as  
\begin{align}
\sum_{k=0}^n\binom{n+1}{k}D_k(t)=t^n. \label{BerId-3} 
\end{align}
Replacing $t$ by $\beta$ we see that $D_k(\beta)$ $(k\ge 0)$ also satisfy the 
recurrent relation that satisfied by $E_k$ $(k\ge 0)$ in 
Eq.\,(\ref{Lma-3.9-peq-4}). Furthermore, by Eqs.\,(\ref{B1}) and (\ref{Def-Dn})
we also have $D_0=1=E_0$. Since every solution of the recurrent relation in 
Eq.\,(\ref{Lma-3.9-peq-4}) is completely determined by the initial 
value for $E_0$, we see that $E_k=D_k(\beta)$ for all $k\ge 1$, whence the lemma follows. 
\epfv 

Note that for the ideal $I$ in Lemma \ref{Dn-Lma} it is easy to see that 
$\delta I$ does not contain any nonzero constant. From this fact 
we immediately have the following 

\begin{corol} \label{Dn-Corol}
Let $I$, $\beta$ be as in Lemma \ref{Dn-Lma}, and $f(x)=\sum_{i=0}^d a_i x^i \in \kx$. Then  $f(x)\in \delta I$, if and only if the following equation holds:
\begin{align} \label{Dn-Corol-eq1}
\sum_{i=0}^d a_i D_i(\beta)c^i=0.
\end{align}
\end{corol}

More generally, for the image under $\delta$ of an ideal $I=u\kx$ with 
$\deg u\ge 2$, we have the following: 

\begin{rmk}\label{Multi-func-descrip}
Assume that $K$ is algebraically closed. Let $u\in \kx$ with 
$d\!:=\deg u\ge 2$, $I=u\kx$, and $r_i$ $(1\le i\le d)$ be all the 
roots of $u$ in $K$. Set $u_{ij}\!:=(x-r_i)(x-r_j)$ and 
$I_{ij}\!:=u_{ij}\kx$ for all $1\le i<j\le d$. Then 
we have 
\begin{align}
I&=\bigcap_{1\le i<j\le d} I_{ij}, \nno\\
\delta I&= \bigcap_{1\le i<j\le d} \delta I_{ij}. \label{Multi-func-descrip-eq1}
\end{align}

On the other hand, let $T_{i}$ ($1\le i\le d-1$) be 
the affine translation of $\kx$ that maps $x$ to $x+r_i$. 
Then $T_i$ commutes with $\delta$ and maps the ideal $I_{ij}$ ($1\le j\le d$) to 
the ideal generated by $x(x-a_{ij})$, where $a_{ij}=r_j-r_i$. Then 
by Corollary \ref{Dn-Corol} and Eq.\,(\ref{Multi-func-descrip-eq1}) above 
we see that the polynomials $f\in \delta I$ up to the translations 
$T_{i}$ ($1\le i\le d-1$) are characterized by a system of 
equations as the one in Eq.\,(\ref{Dn-Corol-eq1}). 
\end{rmk}

Now we are ready to show the following crucial lemma. 
 
\begin{lemma}\label{Quad-Lma}
Let $u(x)=x(x-a)$ for some $a\in K$ and $I=u\kx$.  
Set $\beta\!:=a/c$. Then the following statements hold:
\begin{enumerate}
   \item[$1)$] if $\beta=1$, then  $\delta I=(x)=xK[x]$;  
\item[$2)$] if $\beta=-1$, then  $\delta I=(x+a)\kx$; 
  \item[$3)$] if $\beta=0$, then $\rad(\delta I)=\{0\}$.
\end{enumerate} 
In all the cases above $\delta I$ is a MS of $\kx$.
\end{lemma}

\pf $1)$  By Eq.\,(\ref{BerId-0}) we have 
$B_n(1)=B_n(0)$ for all $n\ge 2$, and by Eq.\,(\ref{Def-Dn}) 
$D_n(1)=0$ for all $n\ge 1$.   
Furthermore by Eq.\,(\ref{Dk-Lma-eq4}), $D_0(1)=1$. 
Then by Lemma \ref{Dn-Lma} the statement follows.

$2)$ Let $I_1$ and $I_2$ be the ideals of $\kx$ generated respectively 
by $x^2-ax$ and $x^2+ax$. Denote by $T$ the affine translation of $\kx$ that 
maps $x$ to $x+a$. Then $T$ maps the principal ideal $I_1$ to 
$I_2$. Since $T$ is a $K$-algebra automorphism of $\kx$ and commutes 
with $\delta$, we have $\delta I_2=\delta (T I_1)=T(\delta I_1)$.
Since $\delta I_1=(x)$ by statement $1)$, we have $\delta I_2=(x+a)$, 
as desired.

Another proof of this statement is to use Eqs.\,(\ref{BerId-5}) 
and (\ref{Def-Dn}) first to show $D_k(-1)=(-1)^k$ 
for all $k\ge 0$, and then apply Lemma \ref{Dn-Lma}.  

$3)$ Assume otherwise and let $0\ne f(x)\in \rad(\delta I)$, i.e., 
$f^m\in \delta I$ for all $m\gg 0$. Since $\delta I$ does not contain 
any nonzero constant, $d:=\deg f(x)\ge 1$. Furthermore,   
we may assume that $f$ is monic,   
and by replacing $f$ by a power of $f$,  
that $f^m\in \delta I$ for all $m\ge 1$.
Write 
\begin{align}
f(x)=x^d+\sum_{i=0}^{d-1}a_i x^i,  \label{Quad-Lma-peq1}
\end{align}
and for all $m\ge 1$, 
\begin{align}
f^m(x)=x^{md}+\sum_{i=0}^{md-1} \Gamma_{m,j} x^i, \label{Quad-Lma-peq2}
\end{align}
where $\Gamma_{m,j}$'s are some polynomials in $a_i$'s over $\bZ$. 

Applying Corollary \ref{Dn-Corol} to $f^m\in \delta I$ $(m\ge 1)$ we get 
\begin{align}\label{MainProp-3.2-peq2}
D_{md}(0)c^{md}+\sum_{j=0}^{md-1} D_j(0) \Gamma_{m,j}  c^j=0. 
\end{align}
Since by Eq.\,(\ref{Dk-Lma-eq4}) we have $D_n(0)=B_n$ for all $n\ge 1$, the equation 
above becomes  
\begin{align} 
B_{md}c^{md}+\sum_{j=0}^{md-1} B_j \Gamma_{m,j} c^j=0. \label{Quad-Lma-peq4}
\end{align}

Next, we make the following reduction.
 Note that for all $m\ge 1$ the equations above are polynomial equations over $\bQ$ 
in $c$ and $a_i$'s. We may apply a similar reduction as 
in the proof of \cite[Theorem $4.1$]{FPYZ} to assume that 
$K$ is a subfield of the algebraic closure $\bar \bQ$ of 
the rational field $\bQ$. Therefore, for each prime $p$ the $p$-valuation 
$\nu_p$ of $\bQ$ can be extended to $K$, 
which we will still denote by $\nu_p$.             
 
Now, by Dirichlet's prime number theorem there exist infinitely many $m\ge 1$ 
such that $md+1$ is a prime number. 
Furthermore, it is well-known in Algebraic Number Theory 
(e.g., see \cite[Theorem 4.1.7]{W}) that for all but finitely primes $p$,  
the values of $\nu_p$ at $c$ and $a_i$ $(0\le i\le d-1)$ are equal to $0$.  
Therefore, we may choose an $m\ge 1$ such that 
the following properties hold: 
{\it 
\begin{enumerate}
  \item[$i)$] $p\!:=md+1$ is an odd prime (In particular, $md$ is even);
 \item[$ii)$] $\nu_p(c)=\nu_p(a_i)=0$ for all $0\le i\le d-1$ such that $a_i\ne 0$.  
\end{enumerate}
}
Consequently, for all $0\le j\le md-1=p-2$ we have $\nu_p(\Gamma_{m,j})\ge 0$ and    
by Corollary \ref{pv4Ber}, $1)$, $\nu_p\big(B_j \big)\ge 0$. 
Then by Eq.\,(\ref{Quad-Lma-peq4}) $\nu_p(B_{p-1}c^{md})\ge 0$. 
Since $\nu_p(c)=0$, we have $\nu_p(B_{p-1})\ge 0$. 
But this contradicts to Corollary \ref{pv4Ber}, $2)$.  
Therefore statement $3)$ follows.
\epfv

\begin{lemma}\label{RepeatRootCases}
Let $I=u\kx$ such that $u$ has at least one repeated root in $\bar K$. 
Then $\rad(\delta I)=\{0\}$. In particular, $\delta I$ is a MS of $\kx$.
\end{lemma}

\pf Let $\bar K$ be the algebraic closure of $K$. We view $\kx$ as a $K$-subalgebra of $\bar K[x]$ in the canonical way and denote by $\bar\delta$ the $\bar K$-linear extension of $\delta$ from $\bar K[x]$ to $\bar K[x]$.

Let $r$ be a repeated root of $u$ in $\bar K$, and $T$ the affine 
translation of $\bar K[x]$ that maps $x$ to $x+r$. Denote by $I_1$ and $I_2$ 
the ideals of $\bar K[x]$ generated respectively by $x^2$ and 
$(x-r)^2$. Then  $TI_1=I_2$. 

Applying Lemma \ref{Quad-Lma}, $3)$ to $\bar \delta$ and 
$I_1$, we see that the radical $\rad(\bar \delta I_1)$ in 
$\bar K[x]$ is equal to $\{0\}$. Since $T$ commutes 
with $\bar \delta$ and is a $\bar K$-algebra 
automorphism, we have 
$\rad(\bar \delta I_2)=\rad\big (T \bar \delta I_1\big)=
T \big(\rad(\bar \delta I_1)\big)=\{0\}$. 
Since $\delta I \subseteq \bar \delta I_2$,  
the radical $\rad(\delta I)$ in $\kx$ is also equal to $\{0\}$, 
whence the lemma follows.
\epfv


Now we are ready to show the last part of 
Theorem \ref{MainThm-2}, i.e., statement $3)$. \\

\underline{\bf Proof of Theorem \ref{MainThm-2},  3)}:
Note first that, if $c=0$, then $\phi=\I$ 
and $\delta=0$, whence the statement 
holds trivially in this case. 
Therefore we assume $c\ne 0$.

Let $u\in \kx$ and $I=(u)$. If $u=0$, then $\delta I=0$, 
whence the theorem holds. If $u\ne 0$ and $\deg u\le 1$, 
then the statement follows from Lemma \ref{Lma-4.1}. 

If $\deg u\ge 2$ and $u$ has a repeated root in the algebraic closure of $K$,  
then by Lemma \ref{RepeatRootCases} we have $\rad(\delta I)=\{0\}$.
Hence $\delta I$ is a MS of $\kx$, and the statement holds.
\epfv

One consequence of Theorem \ref{MainThm-2}, $3)$ is the following 
corollary on the image of the quantum derivation 
$D_h$ (e.g., see \cite{KC}) for all nonzero $h\in K$, 
which is defined by setting for all $f\in \kx$ 
\begin{align}
D_h f(x)\!:=\frac{f(x+h)-f(x)}{h}.
\end{align}

\begin{corol}\label{h-Der}
Let $c$, $\delta$ be fixed as before, $u\in\kx$ and $I=u\kx$. 
Assume that either $u=0$, or $\deg u\le 1$ or $u$ has 
at least one repeated root in the algebraic closure of $K$. 
Then the quantum derivation $D_{h=c}$ 
maps $I$ to a MS of $\kx$.  
\end{corol}

Another remark on the $K$-$\cE$-derivation $\delta$ studied 
in this section is as follows.

Let $S$ be the affine automorphism of $\kx$ that maps 
$x$ to $c^{-1}x$. Then it is easy to check that 
$S^{-1}\delta S=-\Delta$, where $\Delta$ is the so-called 
{\it difference operator} of $\kx$, i.e., $\Delta f=f(x+1)-f(x)$ 
for all $f\in \kx$. Therefore, all the results obtained for $\delta$ 
in this section can also be interpreted as certain results 
on the images of ideals of $\kx$ under 
the difference operator $\Delta$ of $\kx$.

%

\end{document}